\documentclass[12pt,twoside]{article}
\usepackage[mathcal]{eucal}
\usepackage{bbm}
\usepackage[cp1251]{inputenc}   
\usepackage[english]{babel}
\usepackage{amsmath}
\usepackage{amsfonts,amssymb}
\usepackage[letterpaper,hmargin=1in]{geometry}
\usepackage{graphicx}                 

\usepackage{srcltx} 

\newcommand{\av}{v}
\newcommand{\bw}{w}

\newcommand{\rav}{\stackrel{\triangle}{=}}

\newcommand{\mm}[1]{{\mathbb{#1}}}

\newcommand{\epsi}{\varepsilon}
\newcommand{\bo}{\hfill {$\Box$}}

\newtheorem{definition}{Definition}

\newtheorem{theorem}{Theorem}

\begin{document}
\title{One simple remark concerning the uniform value\thanks{Krasovskii Institute of Mathematics and Mechanics,     Russian
Academy of Sciences, 16, S.Kovalevskaja St., 620990, Yekaterinburg, Russia; \ \
Institute of Mathematics and Computer Science, Ural Federal University, 4, Turgeneva St., 620083, Yekaterinburg, Russia}}

\author{Dmitry Khlopin\\
{\it khlopin@imm.uran.ru} 
}
\maketitle

%
\begin{abstract}

The paper is devoted to dynamic games.
We consider a general enough framework, which is not limited to e.g.\ differential games and could accommodate both discrete and continuous time.  Assuming common dynamics, we study two game families  with  total payoffs that
   are defined either as the Ces\`{a}ro average (long run average game family)
 or  Abel average (discounting game family) of the running costs. We study a robust strategy that would provide a near-optimal total payoff for all sufficiently small discounts and for all sufficiently large  planning horizons. Assuming merely the Dynamic Programming Principle,  we prove the following Tauberian theorem: if a  strategy  is  uniformly optimal  for
one of the families
(when discount goes to zero  for discounting games,
when planning horizon goes to infinity in long run average games) and its value functions converge uniformly, then, for the other family, this
strategy is also uniformly optimal and its value functions  converge  uniformly to the
same limit.


\end{abstract}

In dynamic optimization,
most often, the potential infinity of the planning horizon is emulated by considering the running cost averaged with respect to the uniform or exponential distributions (Ces\`{a}ro average or Abel average, respectively). In this paper, we consider the correspondence between the value function of  dynamic games with these  total payoff families. The theorems that describe the connection between the Ces\`{a}ro and Abel averages are long known as Tauberian.

The existence of a limit of the value function when the running cost is averaged with respect to the uniform or exponential distributions means that the value functions are robust with respect to the choice of discount (the planning horizon) as long as it is small (large) enough. In particular, it is exactly this value (asymptotic approach) that is viewed as the game's value for the infinite planning horizon in stochastic statements. Furthermore, in these statements, often enough (see~\cite{BK,MN1981}), one could also find a robust strategy that would provide a near-optimal total payoff for all sufficiently small discounts and for all sufficiently large  planning horizons \cite{SZ,Ziliotto2016,JN non,LaR} (uniform approach).

Within the framework of the asymptotic approach, under mild assumptions, there holds the following Tauberian theorem:
the uniform convergence of the value functions for the running costs averaged with respect to the uniform and/or exponential distributions guarantees,  for the other distribution, the uniform convergence of its value functions to the same limit.
Such a result was proved for stochastic games with finite numbers of states and actions~\cite{MN1981} and for discrete-time optimal control problems~\cite{Lehrer}. It has been relatively recently transferred to the general control problems~\cite{barton}, differential games~\cite{Khlopin3}, and a broad class of stochastic games~\cite{Ziliotto2013}. Then, a Tauberian theorem for all two-person zero-sum games satisfying the Dynamic Programming Principle was proved~\cite{DGAA}.

Surprisingly,  the Tauberian theorem proved in~\cite{DGAA}  lends itself well to the proof of the corresponding Tauberian theorem for the uniform approach. This result is the main contribution of this article.

 {\bf Dynamic system}\
  
  Assume the following items are given:
 \begin{itemize}
   \item   a nonempty set  $\Omega$, the data space;
   \item   a nonempty set $\mm{K}$ of maps from $[0,\infty)$ to $\Omega$;
   \item   a running cost $g\colon\Omega\mapsto [0,1]$.
 \end{itemize}
For each process $z\in \mm{K}$, let the  map $t\mapsto g(z(t))$ be
       be Borel measurable. Now,
  for all positive $\lambda,T>0$,  consider the payoffs
 \begin{eqnarray*}
\av_T(z)&\rav&
 \frac{1}{T}\int_{0}^h g(z(t))\,dt\qquad\forall  z\in\mm{K};\\
\bw_\lambda(z)&\rav&
\lambda\int_{0}^\infty e^{-\lambda t}g(z(t))\,dt\qquad\forall z\in\mm{K}.
\end{eqnarray*}

{\bf On game value maps}\

Denote by $\mathfrak{U}$ the set of all bounded maps from $\Omega$ to $\mm{R}$;  denote  by $\mathfrak{C}$ a non-empty set of maps from $\mm{K}$ to $\mm{R}$.
Hereinafter assume that the set  $\mathfrak{C}$ contains all conceivable payoffs,
and the set  $\mathfrak{U}$ contains all value functions for all games with  payoffs $c\in\mathfrak{C}$.

\begin{subequations}
Let $\mathfrak{C}$  satisfy
the following condition:
\begin{equation}
 Ac+B\in\mathfrak{C}\ \mathrm{  for\ all }\ A\geq 0,B\in\mm{R}\ \mathrm{ if } \ c\in\mathfrak{C}.
 \label{conditions}
\end{equation}
  Hereinafter we assume that $\av_T,\bw_\lambda\in\mathfrak{C}$ for all positive $\lambda,T.$

A map $V$ from $\mathfrak{C}$ to $\mathfrak{U}$ is called  a \emph{game value map} if the following conditions hold:
\begin{eqnarray}
\label{conditions1}
  &\ &V[Ac+B]=A\,V[c]+B\  \textrm{ for all } c\in\mathfrak{C}, A\geq 0, B\in\mm{R},\\
  &\ &V[c_1](\omega)\leq V[c_2](\omega)\ \textrm{ for all }  \omega\in\Omega\ \textrm{ if } c_1(z)\leq c_2(z) \ \textrm{for all }
  z\in\mm{K}.
\label{conditions2}
\end{eqnarray}
\end{subequations}


{\bf On Dynamic Programming Principle}\

 Fix a game value map $V$.
 For all positive $\lambda,T,h>0$ and the game value map $V$, define payoffs $\zeta_{h,T}\colon\mm{K}\to\mm{R}$,
 $\xi_{h,\lambda}\colon\mm{K}\to\mm{R}$ as follows:
 \begin{eqnarray*}
 \zeta_{h,T}(z)&\rav&
 \frac{1}{T+h}\int_{0}^h g(z(t))\,dt+\frac{T}{T+h}{V[\av_T]}(z(h))\qquad\forall  z\in\mm{K};\\
\xi_{h,\lambda}(z)&\rav&
\lambda\int_{0}^h e^{-\lambda t}g(z(t))\,dt+e^{-\lambda h}V[\bw_\lambda](z(h))\qquad\forall z\in\mm{K}.
\end{eqnarray*}

\begin{definition}
Let us say that the payoffs  $\av_T(T>0)$ \textup{(}resp., $\bw_\lambda(\lambda>0)$\textup{)}
	enjoy the \emph{weak Dynamic Programming Principle}   with respect to the  game value map~$V$ if the payoffs
    $\zeta_{h,T}(T,h\in\mathbb{N})$ \textup{(}resp., payoffs $\xi_{h,\lambda}(\lambda>0,h\in\mathbb{N})$\textup{)} lie in~$\mathfrak{C}$ and,
	for every $\varepsilon>0$,
	there exists natural  $N$ such that,
	for all natural $h,T>N$  and positive $\lambda<1/N$,
\begin{eqnarray}
    \big|V[\av_{T+h}](\omega)-V[\zeta_{h,T}](\omega)\big|<\epsi\qquad
	\Big(\big|V[\bw_\lambda](\omega)-V[\xi_{h,\lambda}](\omega)\big|<\epsi\Big)\qquad \forall\omega\in\Omega.\label{lim}
\end{eqnarray}
\end{definition}

In particular, the family of payoffs $\av_T(T>0)$ (resp., $\bw_\lambda(\lambda>0)$)
enjoys {the weak Dynamic Programming Principle} if
$$V[\av_{T+h}]=V[\zeta_{h,T}],\quad\Big(V[\bw_{\lambda}]=V[\xi_{h,\lambda}]\Big)\qquad
\forall h,T\in\mm{N}, \lambda>0.$$

 Note that the technical  requirements $\zeta_{h,T}\in\mathfrak{C}$ and $\xi_{h,\lambda}\in\mathfrak{C}$ can always be provided for by extending $V$ \cite[Lemma~1]{DGAA}. The key requirement in the Dynamic Programming Principle is the uniform approximativity of the value functions $V[\av_{T+h}]$ and $V[\bw_\lambda]$ with $V[\zeta_{h,T}]$ and
 $V[\xi_{h,\lambda}]$, respectively.

 The cornerstone result for this article's main theorem was proved in \cite{DGAA}:
\begin{theorem}\label{normal2}
Let there be given a game value map $V\colon\mathfrak{C}\to\mathfrak{U}$. Let  $\av_T,\bw_\lambda\in\mathfrak{C}$ for all
$\lambda,T>0$.
Assume that payoffs $\av_T (T>0)$ and  payoffs $\bw_\lambda (\lambda>0)$
enjoy the weak Dynamic Programming Principle.

Then, the following two statements are equivalent:
\begin{description}
  \item[$(\imath)$]   The family of functions $V[\av_T]$ $(T>0)$  converges uniformly on $\Omega$ as $T\uparrow\infty$.
  \item[$(\imath\imath)$] The family of functions $V[\bw_\lambda]$ $(\lambda>0)$ converges uniformly on $\Omega$ as $\lambda\downarrow0$.
\end{description}
Moreover, when at least one of these statements hold,  then, for
    both families, the corresponding limits of the value functions
       exist, are uniform in $\omega\in\Omega$, and coincide.
\end{theorem}

Note that in this theorem the requirement of the dynamic programming principle  could not be omitted, see \cite{MathSb}.
The limits must remain uniform unless additional assumptions are made, even in control problems and the stochastic statement;
see the counterexample in e.g.~\cite{barton} and \cite{Ziliotto2013}, respectively.

{\bf On strategies}

Assume that a strategy set $\mathfrak{S}$ is given, and, for each strategy $s\in\mathfrak{S}$, we construct a game value map $V_s\colon\mathfrak{C}\to\mathfrak{U}$.
Consider  the game value map
\begin{equation}
V_{\mathrm{best}}[c](\omega)=\sup_{s\in\mathfrak{S}} V_{s}[c](\omega)\qquad \forall c\in\mathfrak{C},\omega\in\Omega.\label{Vbest}
\end{equation}

\begin{definition}
Let us say that a  strategy $s\in\mathfrak{S}$ is called uniformly optimal for the  payoff family $\av_T({T>0})$ (resp.,~$\bw_\lambda({\lambda>0})$) iff $$  \lim_{T\uparrow \infty}\sup_{\omega\in\Omega}\big|V_{\mathrm{best}}[\av_T]- V_{s}[\av_T]\big|=0\qquad
\Big(\lim_{\lambda\downarrow 0}\sup_{\omega\in\Omega}\big|V_{\mathrm{best}}[\bw_\lambda]- V_{s}[\bw_\lambda]\big|=0\Big).$$
\end{definition}

An unexpected pleasure is the fact that the assumptions of the following theorem  allow the use of any strategy, a strategy of whatever kind: if only it satisfies the Dynamic Programming Principle, we get a uniformly optimal strategy for both payoff families.
\begin{theorem}\label{best}
Let there be given  game value maps $V_s\colon\mathfrak{C}\to\mathfrak{U}$ for each strategy $s\in\mathfrak{S}$, and let the game value map
$V_{\mathrm{best}}$ be defined by the rule \eqref{Vbest}. Let there be given a strategy $s_*\in\mathfrak{S}.$

Assume that the payoffs $\av_T (T>0)$ and  payoffs $\bw_\lambda (\lambda>0)$
enjoy the weak Dynamic Programming Principle with respect to $V_{s_*}$ and $V_{\mathrm{best}}$.
Then, the following conditions are equivalent:
\begin{description}
		\item[$(v)$]
   the strategy $s_*$ is  uniformly  optimal  for the payoff family $\av_T({T>0})$, in addition,
   the functions $V_{\mathrm{best}}[\av_T]$ converge uniformly in $\Omega$ as $T\uparrow \infty$;

  \item[$(w)$]
   the strategy $s_*$ is  uniformly  optimal  for the payoff family $\bw_\lambda({\lambda>0})$, in addition,
   the functions $V_{\mathrm{best}}[\bw_\lambda]$ converge uniformly in $\Omega$ as $\lambda\downarrow 0$;

  \item[$(eq)$]
  all limits in
  \begin{equation}\label{1800}
\lim_{T\uparrow\infty}V_{\mathrm{best}}[\av_T](\omega),\ \lim_{\lambda\downarrow 0}
V_{\mathrm{best}}[\bw_\lambda](\omega),\
\lim_{T\uparrow\infty}V_{s_*}[\av_T](\omega),\ \lim_{\lambda\downarrow 0}
V_{{s_*}}[\bw_\lambda](\omega)
\end{equation}
        exist, are uniform in $\omega\in\Omega$, and coincide.
\end{description}
\end{theorem}

\underline{Proof.}

Note that $(eq)\Rightarrow (v)$ and $(eq)\Rightarrow (w)$ hold by  the definition.
So, it would suffice to prove
  $(v)\Rightarrow (eq)$, $(w)\Rightarrow (eq)$.
 We will prove $(v)\Rightarrow (eq)$; the proof of the last implication $(w)\Rightarrow (eq)$ is  similar.

Let
   the strategy $s_*$ be  uniformly  optimal for  the payoff family $\av_T({T>0})$, and let
   the functions $V_{\mathrm{best}}[\av_T]$ converge uniformly in $\Omega$ as $T\uparrow \infty$ to a function $U_*\in\mathfrak{U}$.
   Then, by the definition of the uniformly  optimal strategy,
     the functions $V_{s_*}[\av_T]$ also converge to $U_*$ uniformly   in $\Omega$ as $T\uparrow \infty$.
    Applying Theorem~\ref{normal2} for the game value map $V_{s_*}$, we see that
    $V_{s_*}[\bw_\lambda]$ converge to $U_*$ uniformly   in $\Omega$ as $\lambda\downarrow 0$.
    Applying this theorem for the game value map $V_{\mathrm{best}}$, we have that
    $V_{\mathrm{best}}[\bw_\lambda]$ converge to $U_*$ uniformly   in $\Omega$ as $\lambda\downarrow 0$.
    Then, in view of the definition  of the uniformly optimal strategy, we see that
    $s_*$ is uniformly optimal for   the payoff family $\bw_\lambda(\lambda>0)$.
 \bo


 \underline{Question~1}.
  Is the uniform  optimality of a strategy  for the payoffs $\av_T(T>0)$ and for the payoffs $\bw_\lambda( \lambda>0)$  still equivalent without the additional requirement on  the uniform convergence of arbitrary value  functions?

  \underline{Question~2}.
  Could the assumption on Dynamic Programming Principle be relaxed?

  \underline{Question~3}. Can there be a topology for $\mathfrak{U}$ that is different from the uniform topology?

  A partial answer to  {Question~2} will be considered for the following important case.

  {\bf The case of one player}\

 For all $\tau\in[0,\infty),z',z''\in\mathbb{K}$ with the property
   $z'(\tau)=z''(0)$ and them only, let us define their concatenation---likewise, a mapping from $[0,\infty)$ to~$\Omega$---as follows:
\begin{equation*}
(z'\diamond_\tau z'') (t)\rav \left\{
 \begin{array} {rcl}        z'(t),       &\mathstrut&    0\leq t\leq\tau;\\
                            z''(t-\tau),
                                         &\mathstrut&     t>\tau.
 \end{array}            \right.\qquad\forall t\geq 0.
\end{equation*}

 Following \cite{barton}, assume that
  \begin{equation*}
  \Gamma(\omega)\rav\{z\in\mathbb{K}\,|\,z(0)=\omega\}\neq\emptyset\qquad \forall \omega\in\Omega; 
  \end{equation*}
  also, assume that the set~$\mathbb{K}$ is closed with respect to concatenation.

  Let $\mathfrak{S}$ be the set of all selectors $\Omega\in\omega\mapsto s[\omega]\in\mathbb{K}$ of the set-valued map $\omega\mapsto\Gamma(\omega)$, i.e., the set of all mappings $\Omega\in\omega\mapsto s[\omega]\in\mathbb{K}$ such that $s[\omega](0)=\omega$ for all $\omega\in\Omega.$

\begin{definition}
Let us say that a strategy $s_*\in\mathfrak{S}$ is {\it stationary-like} (does not change when shifted in time) if
     \begin{equation*}
     s_*[\omega](t+1)=s_*\big[s_*[\omega](1)\big](t)\quad\forall \omega\in\Omega,t\geq 0.
      \end{equation*}
   \end{definition}

   Let $\mathfrak{C}$ and $\mathfrak{U}$ be the sets of all scalar bounded maps from $\mathbb{K}$ and  $\Omega$, respectively.
   For all $s\in\mathfrak{S}$, define $V_s\colon\mathfrak{C}\to\mathfrak{U}$ by the following rule:
   $$V_s[c](\omega)\rav c(s[\omega])\qquad\forall \omega\in\Omega,c\in\mathfrak{C}.$$
   Let us also define $V_{\mathrm{best}}\colon\mathfrak{C}\to\mathfrak{U}$ by \eqref{Vbest}. Then,
   $$V_{\mathrm{best}}[c](\omega)=\sup_{z\in\mathbb{K},z(0)=\omega}c(z)\qquad \forall \omega\in\Omega,c\in\mathfrak{C}.$$
It is   easy to see that $V_{\mathrm{best}}$ and $V_s$, for all $s\in\mathfrak{S}$, are game value maps.

   \begin{theorem}\label{bestbest}
   Assume that  the set $\mathbb{K}$ is closed with respect to concatenation and that
   $\Gamma(\omega)$ is non-empty for all $\omega\in\Omega$.
Let  a strategy $s_*\in\mathfrak{S}$ be stationary-like.

Then, conditions $(v)$, $(w)$, $(eq)$ are equivalent.
\end{theorem}
\underline{Proof.}

It is  easy to see that, since  the strategy $s_*\in\mathfrak{S}$ is  stationary-like, the payoffs $\av_T (T>0)$ and  payoffs $\bw_\lambda (\lambda>0)$
enjoy the weak Dynamic Programming Principle with respect to $V_{s_*}$.

 Since the set $\mathbb{K}$ is closed with respect to concatenation, we have
  \begin{eqnarray*}
  V_{\mathrm{best}}[c](\omega)=\sup_{s\in\mathfrak{S}}c(s[\omega])\geq \sup_{s_0,s_1\in\mathfrak{S}}c(s_0[\omega]\diamond_1 s_1[s_0[\omega](1)]),&\ &\\
  V_{\mathrm{best}}[c](\omega)=\sup_{s\in\mathfrak{S}}c(s[\omega])\geq \sup_{s_0,s_1\in\mathfrak{S}}c(s_0[\omega]\diamond_n s_1[s_0[\omega](n)]),&\ &\forall n\in\mathbb{N}.
  \end{eqnarray*}
 Then,  for all positive $\lambda$ and natural $n,T$, for all $\omega\in\Omega$, we have
  \begin{eqnarray*}
  V_{\mathrm{best}}[\bw_\lambda](\omega)\geq  \sup_{s\in\mathfrak{S}}\Big[ \lambda\int_0^n e^{-\lambda t}g(s[\omega](t))\,dt+
  e^{-\lambda n}V_{\mathrm{best}}[\bw_\lambda](s[\omega](n))\Big],\\
  V_{\mathrm{best}}[\av_{T+n}](\omega)\geq  \sup_{s\in\mathfrak{S}}\Big[ \frac{1}{T+n}\int_0^n g(s[\omega](t))\,dt+\frac{T}{T+n}
  V_{\mathrm{best}}[\av_T](s[\omega](n))\Big].
  \end{eqnarray*}
  Thus,  with respect to the game value map $V_{\mathrm{best}}$, $V_{\mathrm{best}}[\bw_\lambda]$ is a subsolution (see \cite[Definition~1]{DGAA}) for the payoffs $\bw_\lambda (\lambda>0)$, and $V_{\mathrm{best}}[\av_T]$ is a subsolution  for  the payoffs $\av_T (T>0)$.

So, it would again  suffice to prove
  $(v)\Rightarrow (eq)$, $(w)\Rightarrow (eq).$

Let
   the strategy $s_*$ be  uniformly optimal for  the payoff family $\av_T({T>0})$ (resp., $\bw_\lambda (\lambda>0)$), and
let   its value functions $V_{\mathrm{best}}[\av_T]$  converge uniformly in $\Omega$ to a function $U_*\in\mathfrak{U}$.
   Then, by the definition of the uniformly optimal strategy,
the   value functions $V_{s_*}[\av_T]$ (resp., $V_{s_*}[\bw_\lambda]$) also converge uniformly in $\Omega$ to  $U_*\in\mathfrak{U}$.
   Since  the payoffs $\av_T (T>0)$ and  payoffs $\bw_\lambda (\lambda>0)$ enjoy the Dynamic Programming Principle  with respect to $V_{s_*}$, applying Theorem~\ref{normal2} for the game value map $V_{s_*}$, we see that
    $V_{s_*}[\bw_\lambda]$  and $V_{s_*}[\av_T]$ also converge to $U_*$ uniformly   in $\Omega$. So,
   the lower limit of $V_{\mathrm{best}}[\av_T]-U_*$ (resp., of $V_{\mathrm{best}}[\bw_\lambda]-U_*$) is nonnegative.

   On the other hand, applying \cite[Proposition~3]{DGAA} to the subsolution $V_{\mathrm{best}}[\bw_\lambda]$ (applying \cite[Proposition~4]{DGAA} to the subsolution $V_{\mathrm{best}}[\av_T]$)  and for the game value map $V_{\mathrm{best}}$ we find that, for every positive $\varepsilon$, there exists a natural $N$ such that
   $V_{\mathrm{best}}[\bw_\lambda](\omega)\leq V_{\mathrm{best}}[\av_T](\omega)+\varepsilon$ \big(   $V_{\mathrm{best}}[\av_T](\omega)\leq V_{\mathrm{best}}[\bw_\lambda](\omega)+\varepsilon$\big)  for all positive $T>N$, $\lambda=1/N$, and for all $\omega\in\Omega$.
   So,  the upper limit of $U_*-V_{\mathrm{best}}[\av_T]$ (resp., of $U_*-V_{\mathrm{best}}[\bw_\lambda]$) is also nonnegative.

   Thus,     $V_{\mathrm{best}}[\bw_\lambda]$  and $V_{\mathrm{best}}[\av_T]$  converge to $U_*$ uniformly   in $\Omega$.
   Since, see above, $V_{s_*}[\bw_\lambda]$  and $V_{s_*}[\av_T]$ also converge to $U_*$ uniformly   in $\Omega$,
   by the definition, $s_*$ is a uniformly optimal strategy for both the payoffs $\av_T({T>0})$ and $\bw_\lambda (\lambda>0)$.
   \bo

  \end{document}